\documentclass[a4paper,12pt]{article}
\usepackage{amsmath}
\usepackage{amssymb}
\usepackage{enumerate}
\usepackage{color}
 \numberwithin{equation}{section}

\newtheorem{theo}{Theorem}[section]

\newtheorem{lem}{Lemma}[section]

\newtheorem{col}{Corollary}[section]

\newcommand{\be}{\begin{equation}}
\newcommand{\ee}{\end{equation}}
\newcommand\bes{\begin{eqnarray}}
\newcommand\ees{\end{eqnarray}}
\newcommand{\bess}{\begin{eqnarray*}}
\newcommand{\eess}{\end{eqnarray*}}

\begin{document}
\setlength{\baselineskip}{16pt} \pagestyle{myheadings}

\title{The Poincar\'{e} center-focus  problem for  a class of higher order  polynomial differential systems
\thanks{The corresponding author: Zhengxin Zhou, E-mail:zxzhou@yzu.edu.cn }}
\date{\empty}
\author{Zhengxin Zhou\\
{\small  School of Mathematical Sciences, Yangzhou University,
China} } \maketitle

{\bf Abstract.} { %\footnotesize
\small In this paper, I have proved that for a class of polynomial
differential systems of degree $n+1$ ( where $n$ is an arbitrary
positive integer) the composition conjecture is true. I give the
sufficient and necessary conditions for these differential systems
to have a center at origin point by using a different  method from
the previous references.  By this I can obtain all the focal values
of these systems for an arbitrary $n$ and their  expressions  are
 succinct and beautiful. I  believe that the idea and method of this article can be
used to solve the center-focus problem of more high-order polynomial
differential systems.



\section{ Introduction}
Consider the planar analytic differential system
\begin{equation} \label{1.1}
   \left\{\begin{array}{cc}
    x'=-y+p(x,y),\\
    y'= x+q(x,y),
    \end{array}
    \right.
   \end{equation}
  with $p$ and $q$ being polynomials without constant and linear terms,  and
  seek for conditions under which the origin is a center (that is, the critical point at the origin
 is surrounded by closed orbits).
The derivation of conditions for a center is a difficult and
long-standing problem in the theory of nonlinear differential
equations, however due to complexity of the problem
 necessary and sufficient conditions are known only for a very
few families  of polynomial  systems \eqref{1.1}. The  conditions
for a center in the quadratic system have been obtained in
\cite{Kap1,Kap2}, and in \cite{Dul}, \cite{2} and \cite{Ma} the
problem has been solved for systems in which $ p$ and $q$ are cubic
polynomials without quadratic terms. The problem is also solved  for
some familes of cubic systems and systems in the form of the linear
center perturbed by  homogeneous quartic and quintic nonlinearities,
see e.g. \cite{5}, \cite{7},  \cite{8}, \cite{ALS}, \cite{9},
\cite{GR},  \cite{4}, \cite{R-S} and references given there.

 The problem of distinguishing between a center and a focus for polynomial systems
 \eqref{1.1} has an analog for the corresponding
 periodic  differential equations \cite{4, YY}.
To see this we  note that the phase curves of \eqref{1.1} near the
origin (0,0) in polar coordinates $x=r\cos\theta$, $y=r\sin\theta$
are determined by the equation
\begin{equation}\frac{dr}{d\theta}=\frac{\cos\theta p(r\cos\theta,r\sin\theta)+\sin\theta q(r\cos\theta,r\sin\theta)}{1+r^{-1}(\cos\theta q(r\cos\theta,r\sin\theta)-\sin\theta p(r\cos\theta,r\sin\theta))}.\label{1.3}\end{equation}
Therefore,  the planar vector field \eqref{1.1} has a center at $(0,
0)$ if and only if all solutions $r(\theta)$
 of   equation \eqref{1.3} near  the solution $r\equiv 0$ are periodic, $r(0)=r(2\pi).$
In such case it is said that equation \eqref{1.3} has a center at
$r=0$.

If $r^{-1}(\cos\theta q(r\cos\theta,r\sin\theta)-\sin\theta
p(r\cos\theta,r\sin\theta))=f(\theta)\neq -1 $, and
$p=\sum_{i+j=2}^np_{ij}x^iy^j$ and $q=\sum_{i+j=2}^nq_{ij}x^iy^j$,
from \eqref{1.3} we get the polynomial equation
\begin{equation}\frac{dr}{d\theta}=r^2\sum_{i=0}^{n-2}\tilde{A}_i(\theta)r^i,\label{1.4}
\end{equation}
where $\tilde{A}_i(\theta)\ (i=0,1,2,...,n-2)$ are $2\pi$-periodic
functions. Thus, finding
  the conditions for  existence of a center at the origin of system \eqref{1.1} is equivalent to finding   conditions which fulfilment yields $2\pi$-periodicity of all solutions
  of  polynomial equation \eqref{1.4} near  $r=0$ \cite{7}.

 If $p$ and $q$ are homogeneous polynomials of degree $n$, then  the substitution
$$\rho=\frac{r^{n-1}}{1+r^{n-1}(\cos\theta q(\cos\theta,\sin\theta)-\sin\theta p(\cos\theta,\sin\theta))}$$
 transforms equation \eqref{1.3}  into the  Abel equation
\begin{equation}\frac{d\rho}{d\theta}=\rho^2(\tilde{A}_1(\theta)+\tilde{A}_2(\theta)\rho),\label{1.5}\end{equation}
where $\tilde{A}_i(\theta)\ (i=1,2)$ are $2\pi$-periodic functions.
Thus, finding  the center conditions for \eqref{1.1}
 is equivalent to studying   when   Abel  equation \eqref{1.5} has a center at $\rho=0.$ This problem has been investigated   in  \cite{6,9,YY} and some other works.

In \cite{AR} the authors have proved that the rigid system
$$ \left\{\begin{array}{cc}
    x'=-y+x(H_k(x,y)+H_n(x,y)),\\
    y'= x+y(H_k(x,y)+H_n(x,y)),     \end{array}
    \right.$$
with $H_i$ homogeneous polynomial in $x, y$ of degree $i$,  has a
center if and only if  it is reversible for $k=1,2$. But in
\cite{AR} does not prove that this center is a composition
center\cite{6} . In \cite{XL}, the author  has computed by computer
and obtained the center condition for this system with $k, n=2k
(k=2,3,4,5)$, i.e., $n$ is only a finite  number which does not
exceed ten.

 In this paper, I study  the polynomial differential  system of degree $n+1$
\begin{equation} \label{1.6}
   \left\{\begin{array}{cc}
    x'=-y+x(P_1(x,y)+P_n(x,y))=P(x,y),\\
    y'= x+y(P_1(x,y)+P_n(x,y))=Q(x,y),
    \end{array}
    \right.
      \end{equation}
where $P_k(x,y)=\sum_{i+j=k}p_{ij}x^iy^j \,( k=1,n)$,   $p_{i,j} \,
(i,j=0,1,2,..,n)$ are real constants. I give the necessary and
sufficient conditions for the origin of (1.5) to be a center by
using  a different method from the previous references and
\cite{AR}. I only  use a few skillful mathematical analysis
techniques to achieve these results.  Remarkably, the obtained
center conditions contain exactly $[\frac{n}{2}]+1$ relations.   I
apply the obtained results to prove that for  all obtained center
cases the Composition Condition \cite{6}  is valid for the
corresponding periodic differential equation
\begin{equation}\frac{dr}{d\theta}=r(P_1(\cos\theta,\sin\theta)r+P_n(\cos\theta,\sin\theta)r^n). \label{1.7}\end{equation}
This means that  the composition conjecture \cite{6} is true for
this equation. By this I can derive all the focal values  of system
(1.5) with an arbitrary $n$ and the expressions of the focal values
are more succinct and beautiful than the results of computer
calculation.

I am sure that the idea and method of this paper can be used to
solve the center-focus problem of more high-order polynomial
differential systems.

\section{ Main results}

 Alwash and Lloyd \cite{6,7} proved  the following
statement.

\begin{lem}[\cite{6,7}]
If  there exists a differentiable function $u$ of period $2\pi$ such
that
$$\tilde{A}_1(\theta)=u'(\theta)\hat{A}_1(u(\theta)),\,\tilde{A}_2(\theta)=u'(\theta)\hat{A}_2(u(\theta))$$
for some continuous functions $\hat{A}_1$ and $\hat{A}_2$, then the
Abel differential equation
$$\frac{dr}{d\theta}=\tilde{A}_1(\theta)r^2+\tilde{A}_2(\theta)r^3$$
 has a center at $r\equiv 0$.
\end{lem}

The  condition in Lemma 2.1 is called the  {\it{\bf Composition
Condition}}. This is a sufficient but not a necessary condition for
$r=0$ to be a center \cite{5}.

The following statement presents a   generalization of Lemma 2.1.

\begin{lem}[\cite{ZX}]
  If  there exists a differentiable function $u$ of period $2\pi$ such that
$$\tilde{A}_i(\theta)=u'\hat{A}_i(u),\,(i=1,2,...,n)$$
for some continuous functions $\hat{A}_i \,(i=1,2,...,n)$, then the
differential equation
$$\frac{dr}{d\theta}=r\sum_{i=1}^n\tilde{A}_i(\theta)r^i $$
 has a center at $r=0$.
\end{lem}

Consider $2\pi$ -periodic equations:
\begin{equation}
\frac{dr}{d\theta}=r(P_1r+P_{2m+1}r^{2m+1}),\label{2.1}
\end{equation}
\begin{equation}
\frac{dr}{d\theta}=r(P_1r+P_{2m}r^{2m}),\label{2.2}
\end{equation}
where  $P_{2m+1}$ and $P_{2m}$ are polynomials of degree $2m+1$ and
$2m$ with respect to $\cos\theta$ and $\sin\theta$, respectively.

Denote
\begin{equation}P_1=A_1^1\cos\theta+B_1^1\sin\theta,\label{2.3}\end{equation}
 \begin{equation}\bar{P}_1=A_1^1\sin\theta-B_1^1\cos\theta,\label{2.4}\end{equation}
where $A_1^1,\,B_1^1$ are real numbers and ${A_1^1}^2+{B_1^1}^2\neq
0.$

In this paper we will prove the following main results.

\begin{theo} For equation (1.6), $r=0$ is a center if and only if
\begin{equation}
\int_{0}^{2\pi}\bar{P}_1^jP_nd\theta=0,\,(j=0,1,2,...,n).\label{2.5}
\end{equation}
\end{theo}
This statement is equivalent to the following two theorems.

\begin{theo}
For equation (2.1), $r=0$ is a center if and only if
\begin{equation}
\int_{0}^{2\pi}\bar{P}_1^{2j+1}P_{2m+1}d\theta=0,\,(j=0,1,2,...,m).\label{2.6}
\end{equation}
\end{theo}
\begin{theo}
For equation (2.2), $r=0$ is a center if and only if
\begin{equation}
\int_{0}^{2\pi}\bar{P}_1^{2j}P_{2m}d\theta=0,\,(j=0,1,2,...,m).\label{2.7}
\end{equation}
\end{theo}

Before giving the proof of these theorems, we first need the
following lemmas.

\section{Several lemmas}

 In the following, the denote $C_n^k=\frac{n!}{k! (n-k)!}$, $n$ and $k$ are natural numbers and $k\leq n.$

\begin{lem}
Let $$\phi(y)=1+y+y^2+...+y^{n-1}+a_{n+1}y^n+a_{n+2}y^{n+1}+....,$$
$$g(y)=\phi^n(y)=(1+y+y^2+...+y^{n-1}+a_{n+1}y^n+a_{n+2}y^{n+1}+....)^n,$$
where $a_i\, (i=n+1, n+2,...)$ are independent of $y$. Then
\begin{equation}\frac{g^{k}(0)}{k!}=C_{n+k-1}^k=\frac{(n+k-1)!}{k!(n-1)!},
(k=0,1,2,...,n-1).\label{3.1}\end{equation}
\end{lem}

{\bf Proof.} Denote $\phi(y)=\phi_1(y)+\phi_2(y),$ where
$$\phi_1(y)=1+y+y^2+...+y^{n-1},\,\phi_2(y)=a_{n+1}y^n+a_{n+2}y^{n+1}+....,$$then
$$g(y)=(\phi_1(y)+\phi_2(y))^n=\phi_1^n(y)+C_{n}^1\phi_1^{n-1}(y)\phi_2(y)+...+\phi_2^n(y),$$
so it can be seen that the coefficients of $y^k$ $(k=0,1,2,...,n-1)$
of $g(y)$ and  $\phi_1^n(y)$ are the same. On the other hand,
$$\phi_1(y)=1+y+y^2+...+y^{n-1}=\frac{1}{1-y}-(y^n+y^{n+1}+...)=h(y)-\phi_3(y),$$
where $h(y)=\frac{1}{1-y},$ $\phi_3(y)=y^n+y^{n+1}+...\,.$ Then
$$\phi_1^n(y)=h^n(y)-C_n^1h^{n-1}(y)\phi_3(y)+...+(-1)^n\phi_3^n(y),$$
it implies that the coefficients of $y^k$ $(k=0,1,2,...,n-1)$ of
$\phi_1^n(y)$  and  $h^n(y)$ are the same. So,
$$\frac{g^{(k)}(0)}{k!}=\frac{{h^n}^{(k)}(0)}{k!}=\frac{(n+k-1)(n+k-2)\cdot\cdot\cdot(n+1)n}{k!}=C_{n+k-1}^{k},
\, (k=0,1,2,...,(n-1)).$$ Therefore, the relation \eqref{3.1} is
correct. $\Box$

\begin{lem} Let \begin{equation}\lambda_0^0=C_{2m-1}^0, \, \lambda_{k}^{j}=\frac{2}{k-j}\sum_{i=j}^{k-1}\lambda_{i}^j,\,( j=0,1,2,...,k-1),\,\lambda_k^k=C_{2m+k+1}^k-\sum_{j=0}^{k-1}\lambda_{k}^j.
\label{3.2}\end{equation} Then
\begin{equation}\lambda_{k}^j=(k-j+1)C_{2m+j-1}^j,\,(j=0,1,2,...,k)\label{3.3}
\end{equation}
where $k$ is a positive integer. \end{lem}

{\bf Proof.} Taking $k=1,$  by (3.2), we get
$$\lambda_1^0=2\lambda_0^0=2C_{2m-1}^0,\, \lambda_1^1=C_{2m+2}^1-\lambda_1^0=2m=C_{2m}^1,$$
it implies that the  identity (3.3) is correct when $k=1$.

Suppose that the identity (3.3) is true for positive integer  $k$,
 we will show that the identity (3.3) holds for  integer $k+1$.

Indeed, by (3.2) we have
$$\lambda_{k+1}^{j}=\frac{2}{k+1-j}\sum_{i=j}^{k}\lambda_{i}^j $$
$$=\frac{2}{k+1-j}(\lambda_j^j+\lambda_{j+1}^j+...+\lambda_{k}^j)$$
$$=\frac{2}{k+1-j}(1+2+...,+(k-j+1))C_{2m+j-1}^j$$$$=(k-j+2)C_{2m+j-1}^j, \, (j=0,1,2,...,k),$$
and
$$\lambda_{k+1}^{k+1}=C_{2m+k+2}^{k+1}-\sum_{j=0}^{k}\lambda_{k+1}^j$$
$$=C_{2m+k+2}^{k+1}-(\lambda_{k+1}^0+\lambda_{k+1}^1+...+\lambda_{k+1}^k)$$
$$=C_{2m+k+2}^{k+1}-((k+2)C_{2m-1}^0+(k+1)C_{2m}^1+kC_{2m+1}^2+...+3C_{2m+k-2}^{k-1}+2C_{2m+k-1}^k)$$
$$=C_{2m+k+2}^{k+1}-((k+1)C_{2m-1}^0+kC_{2m}^1+(k-1)C_{2m+1}^2+...+2C_{2m+k-2}^{k-1})$$
\begin{equation}-(C_{2m-1}^0+C_{2m}^1+C_{2m+1}^2+...+C_{2m+k-2}^{k-1}+C_{2m+k-1}^k)-C_{2m+k-1}^k.\end{equation}
As, for arbitrary positive integer $n,\,k$ we have \cite{ZM}
$$C_n^k=C_{n-1}^{k-1}+C_{n-1}^k.$$ By this we obtain
\begin{equation}C_{2m+k+2}^{k+1}=C_{2m+k+1}^{k}+C_{2m+k+1}^{k+1}.\end{equation}
By the inductive hypothesis, we get
\begin{equation}C_{2m+k+1}^{k}-((k+1)C_{2m-1}^0+kC_{2m}^1+(k-1)C_{2m+1}^2+...+2C_{2m+k-2}^{k-1})=\lambda_k^k=C_{2m+k-1}^k.\end{equation}
Substituting (3.6) and (3.5) into (3.4) we get
\begin{equation}\lambda_{k+1}^{k+1}=C_{2m+k+1}^{k+1}-(C_{2m-1}^0+C_{2m}^1+C_{2m+1}^2+...+C_{2m+k-2}^{k-1}+C_{2m+k-1}^k).\end{equation}
Similar to (3.5) we have
$$C_{2m+k+1}^{k+1}=C_{2m+k}^{k}+C_{2m+k}^{k+1}=C_{2m+k-1}^{k-1}+C_{2m+k-1}^{k}+C_{2m+k}^{k+1}$$$$=C_{2m+k-2}^{k-2}+C_{2m+k-2}^{k-1}+C_{2m+k-1}^{k}+C_{2m+k}^{k+1}$$
$$=...=C_{2m-1}^0+C_{2m}^1+C_{2m+1}^2+...+C_{2m+k-2}^{k-1}+C_{2m+k-1}^k+C_{2m+k}^{k+1},$$
substituting it into (3.7) we obtain
$$\lambda_{k+1}^{k+1}=C_{2m+k}^{k+1}.$$
Thus, the identity (3.3) holds for  integer $k+1$. By mathematical
induction, the present lemma is correct.$\Box$

\bigskip
\begin{lem}
Let \begin{equation}
P_1\bar{P}_1=A_2\cos2\theta+B_2\sin2\theta,\end{equation}
\begin{equation}\bar{P}_1^{2k}P_1=\sum_{j=0}^k(A_{2k+1}^{2j+1}\cos(2j+1)\theta+B_{2k+1}^{2j+1}\sin(2j+1)\theta),\end{equation}
where
\begin{equation}A_2=-A_{1}^1B_{1}^1,\,B_2=\frac{1}{2}({A_{1}^1}^2-{B_{1}^1}^2),\end{equation}
\begin{equation}A_{2k+1}^1=\frac{2k-1}{2}(-A_2B_{2k-1}^1+B_2A_{2k-1}^1+\frac{1}{3}(-A_2B_{2k-1}^3+B_2A_{2k-1}^3)),\end{equation}
\begin{equation}B_{2k+1}^1=\frac{2k-1}{2}(-A_2A_{2k-1}^1-B_2B_{2k-1}^1+\frac{1}{3}(A_2A_{2k-1}^3+B_2B_{2k-1}^3)),\end{equation}
\begin{equation}A_{2k+1}^{2j+1}=\frac{2k-1}{2}(\frac{1}{2j-1}(-A_2B_{2k-1}^{2j-1}-B_2A_{2k-1}^{2j-1})+\frac{1}{2j+3}(-A_2B_{2k-1}^{2j+3}+B_2A_{2k-1}^{2j+3})),\end{equation}
\begin{equation}B_{2k+1}^{2j+1}=\frac{2k-1}{2}(\frac{1}{2j-1}(A_2A_{2k-1}^{2j-1}-B_2B_{2k-1}^{2j-1})+\frac{1}{2j+3}(A_2A_{2k-1}^{2j+3}+B_2B_{2k-1}^{2j+3})),\end{equation}
$$\,(j=1,2,3,...,k-2),$$
\begin{equation}A_{2k+1}^{2k-1}=\frac{2k-1}{2(2k-3)}(-A_2B_{2k-1}^{2k-3}-B_2A_{2k-1}^{2k-3}),\end{equation}
\begin{equation}B_{2k+1}^{2k-1}=\frac{2k-1}{2(2k-3)}(A_2A_{2k-1}^{2k-3}-B_2B_{2k-1}^{2k-3}),\end{equation}
\begin{equation}A_{2k+1}^{2k+1}=\frac{1}{2}(-A_2B_{2k-1}^{2k-1}-B_2A_{2k-1}^{2k-1}),\end{equation}
\begin{equation}B_{2k+1}^{2k+1}=\frac{1}{2}(A_2A_{2k-1}^{2k-1}-B_2B_{2k-1}^{2k-1}).\end{equation}
Then
\begin{equation}
A_{2k+1}^{2j+1}=\mu_{2k+1}^{2j+1}A_{2j+1}^{2j+1},\label{3.16}\end{equation}
\begin{equation}
B_{2k+1}^{2j+1}=\mu_{2k+1}^{2j+1}B_{2j+1}^{2j+1},\label{3.17}
\end{equation}
where
\begin{equation}
\mu_{2k+1}^{2j+1}=\frac{2k-1}{2}(\frac{2}{2j-1}\mu_{2k-1}^{2j-1}-\frac{\lambda^2}{2(2j+3)}\mu_{2k-1}^{2j+3}),\label{3.18}
\end{equation}
$$\,(j=0,1,2,...,k-1,\,k=1,2,3,...)$$
$\mu_{3}^{1}=\frac{1}{2}\lambda,\,\lambda=\frac{1}{2}({A_{1}^{1}}^2+{B_{1}^{1}}^2),\,\mu_{2k-1}^{-1}=-\mu_{2k-1}^{1}\mu_{3}^{1},\,
\mu_{j}^{j}=1. \, \mu_i^j=0, A_i^j=0,\,B_i^j=0$ for $j>i$, $i,j$ are
nature numbers.
\end{lem}

{\bf Proof.} Now, we  prove identity (3.19) and (3.20) and (3.21)
are true  by mathematical induction.

Taking $k=1$, by (3.10) and (3.11) and (3.12) we get
$$A_3^1=\frac{1}{2}(-A_2B_{1}^1+B_2A_1^1)=\frac{1}{4}({A_{1}^{1}}^2+{B_{1}^{1}}^2)A_1^1=\frac{1}{2}\lambda A_1^1=\mu_3^1A_1^1,$$
$$B_3^1=\frac{1}{2}(-A_2A_{1}^1-B_2B_1^1)=\frac{1}{4}({A_{1}^{1}}^2+{B_{1}^{1}}^2)B_1^1=\frac{1}{2}\lambda B_1^1=\mu_3^1B_1^1,$$
where
$\mu_{3}^{1}=\frac{1}{4}({A_{1}^{1}}^2+{B_{1}^{1}}^2)=\frac{\lambda}{2}.$
This means that the conclusion of the present lemma is true for
$k=1.$

Suppose that the identity (3.19) and (3.20) and (3.21) hold for
positive integer $k$, next we will prove they also are true for
$k+1$, i.e.,
\begin{equation}
A_{2k+3}^{2j+1}=\mu_{2k+3}^{2j+1}A_{2j+1}^{2j+1},\,
(j=0,1,2,...,k)\end{equation}
\begin{equation}B_{2k+3}^{2j+1}=\mu_{2k+3}^{2j+1}B_{2j+1}^{2j+1},(j=0,1,2,...,k)\end{equation}
where
\begin{equation}
\mu_{2k+3}^{2j+1}=\frac{2k+1}{2}(\frac{2}{2j-1}\mu_{2k+1}^{2j-1}-\frac{\lambda^2}{2(2j+3)}\mu_{2k+1}^{2j+3}),\,(j=0,1,2,...,k).
\end{equation}
Indeed, by (3.11) and (3.19) and (3.20) we have
$$A_{2k+3}^1=\frac{2k+1}{2}(-A_2B_{2k+1}^1+B_2A_{2k+1}^1+\frac{1}{3}(-A_2B_{2k+1}^3+B_2A_{2k+1}^3))$$
\begin{equation}=\frac{2k+1}{2}(\mu_{2k+1}^1(-A_2B_1^1+B_2A_1^1)+\frac{1}{3}\mu_{2k+1}^3(-A_2B_3^3+B_2A_3^3)).\end{equation}
Using (3.17) and (3.18) we have
\begin{equation}
A_3^3=\frac{1}{2}(-A_2B_1^1-B_2A_1^1),\,\,
B_3^3=\frac{1}{2}(A_2A_1^1-B_2B_1^1).
\end{equation}
Substituting (3.26) into (3.25) we get
$$A_{2k+3}^1=\frac{2k+1}{2}(2\mu_{2k+1}^1\mu_3^1-\frac{\lambda^2}{6}\mu_{2k+1}^3)A_1^1=\mu_{2k+3}^1A_1^1.$$
Similarly, we  obtain
$$B_{2k+3}^1=\frac{2k+1}{2}(2\mu_{2k+1}^1\mu_3^1-\frac{\lambda^2}{6}\mu_{2k+1}^3)B_1^1=\mu_{2k+3}^1B_1^1,$$
where
$$\mu_{2k+3}^1=\frac{2k+1}{2}(2\mu_{2k+1}^1\mu_3^1-\frac{\lambda^2}{6}\mu_{2k+1}^3)=\frac{2k+1}{2}(-2\mu_{2k+1}^{-1}-\frac{\lambda^2}{6}\mu_{2k+1}^3),$$
i.e., the identities (3.22) and (3.23) and (3.24) are true for
$j=0$.

For $0<j\leq k-2$, by (3.13) and (3.19) and (3.20) we have
$$A_{2k+3}^{2j+1}=\frac{2k+1}{2}(\frac{1}{2j-1}(-A_2B_{2k+1}^{2j-1}-B_2A_{2k+1}^{2j-1})+\frac{1}{2j+3}(-A_2B_{2k+1}^{2j+3}+B_2A_{2k+1}^{2j+3})),$$
\begin{equation}=\frac{2k+1}{2}(\frac{1}{2j-1}\mu_{2k+1}^{2j-1}(-A_2B_{2j-1}^{2j-1}-B_2A_{2j-1}^{2j-1})+\frac{1}{2j+3}\mu_{2k+1}^{2j+3}(-A_2B_{2j+3}^{2j+3}+B_2A_{2j+3}^{2j+3})).\end{equation}
Using (3.17) and (3.18) we have
$$A_{2j+1}^{2j+1}=\frac{1}{2}(-A_2B_{2j-1}^{2j-1}-B_2A_{2j-1}^{2j-1}),$$
$$A_{2j+3}^{2j+3}=\frac{1}{2}(-A_2B_{2j+1}^{2j+1}-B_2A_{2j+1}^{2j+1}),$$
$$B_{2j+3}^{2j+3}=\frac{1}{2}(A_2A_{2j+1}^{2j+1}-B_2B_{2j+1}^{2j+1}),$$
substituting these relations into (3.27) we obtain
$$A_{2k+3}^{2j+1}=\frac{2k+1}{2}(\frac{2}{2j-1}\mu_{2k+1}^{2j-1}A_{2j+1}^{2j+1}-\frac{\lambda^2}{2(2j+3)}\mu_{2k+1}^{2j+3}A_{2j+1}^{2j+1})=\mu_{2k+3}^{2j+1}A_{2j+1}^{2j+1},$$
$$\mu_{2k+3}^{2j+1}=\frac{2k+1}{2}(\frac{2}{2j-1}\mu_{2k+1}^{2j-1}-\frac{\lambda^2}{2(2j+3)}\mu_{2k+1}^{2j+3}).$$
Similar as above and using (3.14) and (3.17)--(3.20) and
$$B_{2j+1}^{2j+1}=\frac{1}{2}(A_2A_{2j-1}^{2j-1}-B_2B_{2j-1}^{2j-1}),$$
we get
$$B_{2k+3}^{2j+1}=\frac{2k+1}{2}(\frac{2}{2j-1}\mu_{2k+1}^{2j-1}B_{2j+1}^{2j+1}-\frac{\lambda^2}{2(2j+3)}\mu_{2k+1}^{2j+3}B_{2j+1}^{2j+1})=\mu_{2k+3}^{2j+1}B_{2j+1}^{2j+1}.$$
So, the identities (3.22) and (3.23) and (3.24) are held for $0\leq
j\leq k-2.$

Using (3.17) and (3.18) we get
$$A_{2k-1}^{2k-1}=\frac{1}{2}(-A_2B_{2k-3}^{2k-3}-B_2A_{2k-3}^{2k-3}),\,B_{2k-1}^{2k-1}=\frac{1}{2}(A_2A_{2k-3}^{2k-3}-B_2B_{2k-3}^{2k-3}),$$
$$-A_2B_{2k+1}^{2k+1}+B_2A_{2k+1}^{2k+1}=-\frac{1}{2}(A_2^2+B_2^2)A_{2k-1}^{2k-1}=-\frac{\lambda^2}{2}A_{2k-1}^{2k-1},$$
$$A_2A_{2k+1}^{2k+1}+B_2B_{2k+1}^{2k+1}=-\frac{1}{2}(A_2^2+B_2^2)B_{2k-1}^{2k-1}=-\frac{\lambda^2}{2}B_{2k-1}^{2k-1}.$$
Applying these relations and (3.13) and (3.19) and (3.20) we get
$$A_{2k+3}^{2k-1}=\frac{2k+1}{2}(\frac{1}{2k-3}(-A_2B_{2k+1}^{2k-3}-B_2A_{2k+1}^{2k-3})+\frac{1}{2k+1}(-A_2B_{2k+1}^{2k+1}+B_2A_{2k+1}^{2k+1}))$$
$$=\frac{2k+1}{2}(\frac{1}{2k-3}\mu_{2k+1}^{2k-3}(-A_2B_{2k-3}^{2k-3}-B_2A_{2k-3}^{2k-3})-\frac{\lambda^2}{2(2k+1)}A_{2k-1}^{2k-1})$$
$$=\frac{2k+1}{2}(\frac{2}{2k-3}\mu_{2k+1}^{2k-3}A_{2k-1}^{2k-1}-\frac{\lambda^2}{2(2k+1)}A_{2k-1}^{2k-1})=\mu_{2k+3}^{2k-1}A_{2k-1}^{2k-1},$$
and
$$B_{2k+3}^{2k-1}=\frac{2k+1}{2}(\frac{1}{2k-3}(A_2A_{2k+1}^{2k-3}-B_2B_{2k+1}^{2k-3})+\frac{1}{2k+1}(A_2A_{2k+1}^{2k+1}+B_2B_{2k+1}^{2k+1}))$$
$$=\frac{2k+1}{2}(\frac{1}{2k-3}\mu_{2k+1}^{2k-3}(A_2A_{2k-3}^{2k-3}-B_2B_{2k-3}^{2k-3})-\frac{\lambda^2}{2(2k+1)}B_{2k-1}^{2k-1})$$
$$=\frac{2k+1}{2}(\frac{2}{2k-3}\mu_{2k+1}^{2k-3}B_{2k-1}^{2k-1}-\frac{\lambda^2}{2(2k+1)}B_{2k-1}^{2k-1})=\mu_{2k+3}^{2k-1}B_{2k-1}^{2k-1},$$
where
$$\mu_{2k+3}^{2k-1}=\frac{2k+1}{2}(\frac{2}{2k-3}\mu_{2k+1}^{2k-3}-\frac{\lambda^2}{2(2k+1)}).$$
Thus, the identities (3.22) and (3.23) and (3.24) are held for $ j=
k-1.$

By (3.15)$--$(3.20) we have
$$A_{2k+3}^{2k+1}=\frac{2k+1}{2(2k-1)}(-A_2B_{2k+1}^{2k-1}-B_2A_{2k+1}^{2k-1})$$
$$=\frac{2k+1}{2(2k-1)}\mu_{2k+1}^{2k-1}(-A_2B_{2k-1}^{2k-1}-B_2A_{2k-1}^{2k-1})$$
$$=\frac{2k+1}{2k-1}\mu_{2k+1}^{2k-1}A_{2k+1}^{2k+1}=\mu_{2k+3}^{2k+1}A_{2k+1}^{2k+1},$$
and
$$B_{2k+3}^{2k+1}=\frac{2k+1}{2k-1}\mu_{2k+1}^{2k-1}B_{2k+1}^{2k+1}=\mu_{2k+3}^{2k+1}B_{2k+1}^{2k+1},$$
where $\mu_{2k+3}^{2k+1}=\frac{2k+1}{2k-1}\mu_{2k+1}^{2k-1}.$ Thus,
the relations (3.22) (3.23) (3.24) are true for $j=k.$

Therefore, the identity (3.19) and (3.20) and (3.21) hold for
positive integer $k+1$.

In summary,  by mathematical induction, the identities (3.19) and
(3.20) and (3.21) are correct for $j=0,1,2,...,k-1,\,k=1,2,3,...$ .
$\Box$

\begin{lem}
If \begin{equation}
\int_{0}^{2\pi}\bar{P}_1^{2j+1}P_{2m+1}d\theta=0,\,(j=0,1,2,...,m),\label{3.28}
\end{equation}
then
\begin{equation}
P_{2m+1}=P_1(\eta_0+\eta_2\bar{P}_1^2+...+\eta_{2m}\bar{P}_1^{2m}),\label{3.29}
\end{equation}
where $\eta_{2j}\,(j=0,1,2,...,m)$ are constants.\end{lem}

{\bf Proof.} As $P_{2m+1}$ is a polynomial of degree $2m+1$ with
respect to $\cos\theta$ and $\sin\theta$, its Fourier expansion can
be expressed as
$$P_{2m+1}=\hat{P}_1+\hat{P}_3+...+\hat{P}_{2m+1},$$
where
$$\hat{P}_{2j+1}=a_{2j+1}\cos(2j+1)\theta+b_{2j+1}\sin(2j+1)\theta,$$
$$a_{2j+1}=\frac{1}{\pi}\int_{0}^{2\pi}\cos(2j+1)\theta P_{2m+1}d\theta,\,b_{2j+1}=\frac{1}{\pi}\int_{0}^{2\pi}\sin(2j+1)\theta P_{2m+1}d\theta,\,(j=0,1,2,...,m).$$

 First of all,  we prove that
\begin{equation}
\hat{P}_{2j+1}=P_1(\eta_{2j+1}^0+\eta_{2j+1}^2\bar{P}_1^2+...+\eta_{2j+1}^{2j}\bar{P}_1^{2j}),\,
(j=0,1,2,...,m)\label{3.30}
\end{equation}
and
\begin{equation}
A_{2j+1}^{2i+1}b_{2i+1}-B_{2j+1}^{2i+1}a_{2i+1}=0,\,(i=0,1,2,...,j,\,
j=0,1,2,...m),\label{3.31}
\end{equation}
where the denotes $A_i^j,\, B_i^j$ are the same as in Lemma 3.3,
$\eta_{2j+1}^{2i} \, (i=0,1,2,...,j,\, j=0,1,2,...m)$ are constants.

   Taking $j=0$, from (3.28) we get
 $$\int_{0}^{2\pi}\bar{P}_1P_{2m+1}d\theta=0,$$
substituting (2.4) into it we obtain
$$\int_{0}^{2\pi}(A_{1}^1 \sin\theta P_{2m+1}-B_1^1\cos\theta P_{2m+1})d\theta=0,$$
it follows
\begin{equation}
A_1^1b_1-B_1^1a_1=0.\label{3.32}
\end{equation}

{\bf Case 1.} If $A_1^1\neq 0$, from (3.32) we get
$b_1=\frac{a_1}{A_1^1}B_1^1$ and
$$\hat{P}_1=a_1\cos\theta+b_1\sin\theta=\frac{a_1}{A_1^1}P_1=\eta_1^0 P_1, $$
where $\eta_1^0=\frac{a_1}{A_1^1}.$

{\bf Case 2.} If $A_1^1= 0$,  $B_1^1\neq 0,$ then
$P_1=B_1^1\sin\theta.$ By (3.32) we get $a_1=0$ and
$$\hat{P}_1=b_1\sin\theta=\frac{b_1}{B_1^1}P_1=\eta_1^0 P_1, $$
where $\eta_1^0=\frac{b_1}{B_1^1}$.

Thus, the relation (3.30) and (3.31) hold for $j=0$.

Suppose that the relations (3.30)  and (3.31) are true for $j=k-1$,
i.e.,
\begin{equation}
\hat{P}_{2k-1}=P_1(\eta_{2k-1}^0+\eta_{2k-1}^2\bar{P}_1^2+...+\eta_{2k-1}^{2k-2}\bar{P}_1^{2k-2})\label{3.33}
\end{equation}
and \begin{equation}
A_{2k-1}^{2i-1}b_{2i-1}-B_{2k-1}^{2i-1}a_{2i-1}=0,\,(i=1,2,...,k).\label{3.34}
\end{equation}
Next, we will check that these relations for $j=k$  hold, too.

Indeed, using (3.9) we get
$$\bar{P}_1^{2k+1}=(2k+1)\int\bar{P}_1^{2k}P_1d\theta$$
$$=(2k+1)\sum_{j=0}^k(\frac{1}{2j+1}A_{2k+1}^{2j+1}\sin(2j+1)\theta-\frac{1}{2j+1}B_{2k+1}^{2j+1}\cos(2j+1)\theta).$$
So, from
$$\int_{0}^{2\pi}\bar{P}_1^{2k+1}P_{2m+1}d\theta=0$$
follows
\begin{equation}A_{2k+1}^1b_1-B_{2k+1}^1a_1+\frac{1}{3}(A_{2k+1}^3b_3-B_{2k+1}^3a_3)+...+\frac{1}{2k+1}(A_{2k+1}^{2k+1}b_{2k+1}-a_{2k+1}B_{2k+1}^{2k+1})=0.\label{3.35}\end{equation}
By lemma 3.3 and (3.34) we have
\begin{equation}A_{2k+1}^{2j+1}b_{2j+1}-B_{2k+1}^{2j+1}a_{2j+1}=\mu_{2k+1}^{2j+1}(A_{2j+1}^{2j+1}b_{2j+1}-B_{2j+1}^{2j+1}a_{2j+1})=0,\,(j=0,1,2,...,k-1).\end{equation}
Substituting (3.36) into (3.35) we obtain
\begin{equation}A_{2k+1}^{2k+1}b_{2k+1}-B_{2k+1}^{2k+1}a_{2k+1}=0.\end{equation}
Thus, the relation (3.31) is held for $j=k$. Next, we will check
that the relation (3.30) is true for $j=k.$

{\bf Case $1^0$}. If $A_{2k+1}^{2k+1}\neq 0,$  using (3.37)  we get
$$b_{2k+1}=\frac{a_{2k+1}}{A_{2k+1}^{2k+1}}B_{2k+1}^{2k+1}.$$ By
(3.9) and Lemma 3.3 and (3.36)  we get
$$\hat{P}_{2k+1}=a_{2k+1}\cos(2k+1)\theta+b_{2k+1}\sin(2k+1)\theta$$$$=\frac{a_{2k+1}}{A_{2k+1}^{2k+1}}(A_{2k+1}^{2k+1}\cos(2k+1)\theta+B_{2k+1}^{2k+1}\sin(2k+1)\theta)$$
$$=\frac{a_{2k+1}}{A_{2k+1}^{2k+1}}(\bar{P}_1^{2k}P_1-\sum_{j=0}^{k-1}(A_{2k+1}^{2j+1}\cos(2j+1)\theta+B_{2k+1}^{2j+1}\sin(2j+1)\theta))$$
\begin{equation}=\frac{a_{2k+1}}{A_{2k+1}^{2k+1}}(\bar{P}_1^{2k}P_1-\sum_{j=0}^{k-1}\mu_{2k+1}^{2j+1}(A_{2j+1}^{2j+1}\cos(2j+1)\theta+B_{2j+1}^{2j+1}\sin(2j+1)\theta)).\label{3.38}\end{equation}
By (2.3) and (3.9) and Lemma 3.3 we have
$$P_1=A_1^1\cos\theta+B_1^1\sin\theta,$$
$$\bar{P}_1^2P_1=A_3^1\cos\theta+B_3^1\sin\theta+A_3^3\cos3\theta+B_3^3\sin3\theta=\mu_3^1(A_1^1cos\theta+B_1^1\sin\theta)+A_3^3\cos3\theta+B_3^3\sin3\theta,$$
so
$$A_3^3\cos3\theta+B_3^3\sin3\theta=\bar{P}_1^2P_1-\mu_3^1P_1=P_1(\alpha_3^0+\bar{P}_1^2),$$
where $\alpha_3^0=-\mu_3^1.$ Suppose that
\begin{equation}A_{2k-3}^{2k-3}\cos(2k-3)\theta+B_{2k-3}^{2k-3}\sin(2k-3)\theta=P_1(\bar{P}_1^{2k-4}+\alpha_{2k-3}^{2k-6}\bar{P}_1^{2k-6}+...+\alpha_{2k-3}^2\bar{P}_1^2+\alpha_{2k-3}^0),\label{3.39}\end{equation}
where $\alpha_{2k-3}^{2j} \, (j=0,1,2,..,k-3)$ are constants. Since
$$\bar{P}_1^{2k-2}P_1=A_{2k-1}^1\cos\theta+B_{2k-1}^1\sin\theta+A_{2k-1}^3\cos3\theta+B_{2k-1}^3\sin3\theta+...$$$$+A_{2k-1}^{2k-1}\cos(2k-1)\theta+B_{2k-1}^{2k-1}\sin(2k-1)\theta$$
$$=\mu_{2k-1}^1(A_1^1cos\theta+B_1^1\sin\theta)+\mu_{2k-1}^3(A_3^3\cos3\theta+B_3^3\sin3\theta)+...$$$$
+\mu_{2k-1}^{2k-3}(A_{2k-3}^{2k-3}\cos(2k-3)\theta+B_{2k-3}^{2k-3}\sin(2k-3)\theta)+A_{2k-1}^{2k-1}\cos(2k-1)\theta+B_{2k-1}^{2k-1}\sin(2k-1)\theta,$$
 by (3.39) we get
$$A_{2k-1}^{2k-1}\cos(2k-1)\theta+B_{2k-1}^{2k-1}\sin(2k-1)\theta=\bar{P}_1^{2k-2}P_1-(\mu_{2k-1}^1(A_1^1cos\theta+B_1^1\sin\theta)$$$$+\mu_{2k-1}^3(A_3^3\cos3\theta+B_3^3\sin3\theta)+...
+\mu_{2k-1}^{2k-3}(A_{2k-3}^{2k-3}\cos(2k-3)\theta+B_{2k-3}^{2k-3}\sin(2k-3)))$$
$$=\bar{P}_1^{2k-2}P_1-(\mu_{2k-1}^1P_1+\mu_{2k-1}^3P_1(\bar{P}_1^2+\alpha_3^0)+...+\mu_{2k-1}^{2k-3}P_1(\bar{P}_1^{2k-4}+\alpha_{2k-3}^{2k-6}\bar{P}_1^{2k-6}+...+\alpha_{2k-3}^2\bar{P}_1^2+\alpha_{2k-3}^0))$$
$$=P_1(\bar{P}_1^{2k-2}+\alpha_{2k-1}^{2k-4}\bar{P}_1^{2k-4}+...+\alpha_{2k-1}^2\bar{P}_1^2+\alpha_{2k-1}^0),$$
where $\alpha_{2k-1}^{2j} \, (j=0,1,2,..,k-2)$ are constants.  Thus,
the relation (2.39) holds for $k+1$. By mathematical induction, the
relation (2.39) is correct for any integer $k$.

Substituting (3.39) into (3.38) we obtain

\begin{equation}\hat{P}_{2k+1}=P_1\sum_{j=0}^k\eta_{2k+1}^{2j}\bar{P}_1^{2j},\end{equation}
where $\eta_{2k+1}^{2j} \, (j=0,1,2,..,k)$ are constants.

{\bf Case $2^0$.} If $A_{2k+1}^{2k+1}=0,\, B_{2k+1}^{2k+1}\neq 0,$
by (3.37) we get $a_{2k+1}=0,$ similar as above we get
$$\hat{P}_{2k+1}=b_{2k+1}\sin(2k+1)\theta=\frac{b_{2k+1}}{B_{2k+1}^{2k+1}}B_{2k+1}^{2k+1}\sin(2k+1)\theta$$
$$=\frac{b_{2k+1}}{B_{2k+1}^{2k+1}}(\bar{P}_1^{2k}P_1-\sum_{j=0}^{k-1}(A_{2k+1}^{2j+1}\cos(2j+1)\theta+B_{2k+1}^{2j+1}\sin(2j+1)\theta))$$
$$=P_1\sum_{j=0}^k\eta_{2k+1}^{2j}\bar{P}_1^{2j},$$
where $\eta_{2k+1}^{2j} \, (j=0,1,2,..,k)$ are constants.

{\bf Case $3^0$.} If $A_{2k+1}^{2k+1}=0,\, B_{2k+1}^{2k+1}= 0,$ then
$\hat{P}_{2k+1}=0,$ taking $\eta_{2k+1}^{2j}=0 \, (j=0,1,2,..,k)$,
then the relation (3.40) is true, too.

Thus the relation (3.30) holds for $j=k.$ By mathematical induction,
the relations (3.30) and (3.31) are true for arbitrary integer $j$.
Therefore, under the  hypothetical condition (3.28) we have
 $$P_{2m+1}=\hat{P}_1+\hat{P}_3+...+\hat{P}_{2m+1}$$
$$=P_1\sum_{k=0}^m\sum_{j=0}^{k}\eta_{2k+1}^{2j}\bar{P}_1^{2j}$$
$$=P_1(\eta_1^0+(\eta_3^0+\eta_3^2\bar{P}_1^2)+...+(\eta_{2m+1}^0+\eta_{2m+1}^2\bar{P}_1^2+...+\eta_{2m+1}^{2m}\bar{P}_1^{2m}))$$
$$=P_1\sum_{k=0}^m\eta_{2k}\bar{P}_1^{2k}.$$
where $\eta_{2k}\,(k=0,1,2...,m)$ are constants. Therefore, the
conclusion of the present lemma is correct. $\Box$

\bigskip

\section{Proof of main results}

Now we prove that the conclusion of Theorem 2.2 is correct.

{\bf Proof.} {\bf Necessity:}

Let $ r(\theta, c)$  be the solution of (2.1) such that $ r(0,c)=c
\,(0<c\ll 1).$  We write
$$
r(\theta,c)=\sum_{n=1}^{\infty}a_n(\theta)c^n,
$$
where $a_1(0)=1$ and $a_n(0)=0$ for $n>1$. The origin of (2.1) is a
center if and only if $r(\theta+2\pi,c)=r(\theta,c)$, i.e.,
$a_1(2\pi)=1,\,a_n(2\pi)=0\,(n=2,3,4, ...)$ \cite{6,7}.

 Substituting $r(\theta,c)$  into (2.1) we obtain

\begin{equation}\sum_{n=1}^{\infty}a_n'(\theta)c^n=P_1(\theta)(\sum_{n=1}^{\infty}a_n(\theta)c^n)^2+P_{2m+1}(\theta)(\sum_{n=1}^{\infty}a_n(\theta)c^n)^{2m+2}.\label{4.1}\end{equation}

Equating the corresponding coefficients of $c^n$ of (4.1) yields
  $$a_1'(\theta)=0, \,a_1(0)=1,$$
\begin{equation}a_{k}'(\theta)=P_1\sum_{i+j=k}a_ia_j,\,a_k(0)=0,\,(k=2,3,...,2m+1).$$
Solving these $2m+1$ equations we get
$$a_k(\theta)=\tilde{P}_1^{k-1},\,\,\,
(k=1,2,3,...,2m+1),\end{equation} where
$$\tilde{P}_1=\int_0^{\theta}P_1d\theta=\bar{P}_1+B_1^1,\,\bar{P}_1=A_1^1\sin\theta-B_1^1\cos\theta.$$
Obviously,  $a_k(\theta+2\pi)=a_k(\theta)\,(k=1,2,...,2m+1).$  Using
(4.2) and  taking $y=\tilde{P}_1r$ in  Lemma 3.1, it can be seen
that the coefficient of $c^k$ of
$(\sum_{n=1}^{\infty}a_n(\theta)c^{n-1})^{2m+2}$ is
$C_{2m+k+1}^{k}\tilde{P}_1^k,$  by (4.1) we get
\begin{equation}a'_{2m+2+k}=P_1\sum_{i+j=2m+2+k}a_ia_j+P_{2m+1}C_{2m+k+1}^k\tilde{P}_1^k, \, a_{2m+2+k}(0)=0,\,(k=0,1,2,...,2m+1).\label{4.3}\end{equation}
 Taking $k=0$ and solving equation (4.3)  we get
\begin{equation}a_{2m+2}=\tilde{P}_1^{2m+1}+\lambda_0^0\bar{P}_{2m+1},\end{equation}
 where $\lambda_0^0=C_{2m+1}^0=C_{2m-1}^0=1,\, \bar{P}_{2m+1}=\int_0^{\theta}P_{2m+1}d\theta.$

Suppose that, for  arbitrary integer $k$ $\,(k<2m+1)$, solving
equation (4.3) we get
\begin{equation}
a_{2m+2+k}=\tilde{P}_1^{2m+k+1}+\sum_{j=0}^k\lambda_k^j\tilde{P}_1^{k-j}\overline{\tilde{P}_1^jP_{2m+1}},\label{4.5}
\end{equation}
where $\lambda_k^j\,(j=0,1,2,...,k)$ are the same as in Lemma 3.2,
$\overline{\tilde{P}_1^jP_{2m+1}}=\int_0^{\theta}\tilde{P}_1^jP_{2m+1}d\theta.$

In the following we will check that (4.5) is true for integer $k+1$.

Indeed, by $(4.2) -- (4.5)$, we have
$$a'_{2m+k+3}=P_1\sum_{i+j=2m+k+3}a_ia_j+P_{2m+1}C_{2m+k+2}^{k+1}\tilde{P}_1^{k+1}$$
$$=(2m+k+2)P_1\tilde{P}_1^{2m+k+1}+2P_1\sum_{i=0}^k\lambda_k^i\tilde{P}_1^{k-i}\overline{\tilde{P}_1^kP_{2m+1}}+$$
$$+2P_1\tilde{P}_1\sum_{i=0}^{k-1}\lambda_{k-1}^i\tilde{P}_1^{k-1-i}\overline{\tilde{P}_1^{k-1}P_{2m+1}}+...\,\,\,...\,\,...+$$
$$+2P_1\tilde{P}_1^{k-1}(\lambda_1^0\tilde{P}_1\bar{P}_{2m+1}+\lambda_1^1\overline{\tilde{P}_1P_{2m+1}})+2P_1\tilde{P}_1^{k}\lambda_0^0\bar{P}_{2m+1}+C_{2m+k+2}^{k+1}P_{2m+1}\tilde{P}_1^{k+1}$$
$$=(\tilde{P}_1^{2m+k+2})'+2\sum_{j=0}^k\sum_{i=j}^k\lambda_{i}^jP_1\tilde{P}_1^{k-j}\overline{\tilde{P}^j_1P_{2m+1}}+C_{2m+k+2}^{k+1}P_{2m+1}\tilde{P}_1^{k+1}$$
$$=(\tilde{P}_1^{2m+k+2})'+\sum_{j=0}^k\sum_{i=j}^k\frac{2}{k-j+1}\lambda_{i}^j(\tilde{P}_1^{k-j+1}\overline{\tilde{P}^j_1P_{2m+1}})'+(C_{2m+k+2}^{k+1}-\sum_{j=0}^k\lambda_{k+1}^j)P_{2m+1}\tilde{P}_1^{k+1}.$$
By Lemma 3.2, we obtain
$$a_{2m+k+3}=\tilde{P}_1^{2m+k+2}+\sum_{j=0}^{k+1}\lambda_{k+1}^j\tilde{P}_1^{k-j+1}\overline{\tilde{P}^j_1P_{2m+1}},$$
where $\lambda_{k+1}^j\,(j=0,1,2,...,k+1)$ are the same as in Lemma
3.2. By mathematical induction,  the identity (4.5) is correct for
$k=0,1,2,...,2m+1.$

 As $P_1$ and $P_{2m+1}$ are odd degree polynomial with respect to  $\cos\theta,\, \sin\theta$, so $\tilde{P}_1^k, \,\tilde{P}_{2m+1}^k,$
and $\bar{P}_1^k$ and
$\overline{\tilde{P}_1^{2k}P_{2m+1}}=\int_{0}^{\theta}\tilde{P}_1^{2k}P_{2m+1}d\theta$
are $2\pi$-periodic functions.  Thus, from $a_{2m+2+k}(2\pi)=0$
yields
\begin{equation}\int_{0}^{2\pi}\tilde{P}_1^{j}P_{2m+1}d\theta=0,\,
(j=0,1,2,..., k).\end{equation} Since,
$\tilde{P}_1=\bar{P}_1+B_1^1,$ $B_1^1$ is a constant,
$$\tilde{P}_1^j=(\bar{P}_1+B_1^1)^j=\bar{P}_1^j+C_j^1B_1^1\bar{P}^{j-1}+C_{j}^2{B_1^1}^2\bar{P}_1^{j-2}+...+{B_1^1}^j,$$
so (4.6) is equivalent to
\begin{equation}
\int_{0}^{2\pi}\bar{P}_1^{j}P_{2m+1}d\theta=0,\,(j=0,1,2,...k).
\end{equation}
As $\bar{P}_1$ and $P_{2m+1}$ are odd polynomials with respect to
$\cos\theta$ and $\sin\theta$,
$\int_{0}^{2\pi}\bar{P}_1^{2i}P_{2m+1}d\theta=0.$ Thus, from
relation (4.7) implies
\begin{equation}
\int_{0}^{2\pi}\bar{P}_1^{2k+1}P_{2m+1}d\theta=0,\,(k=0,1,2,...,m).
\end{equation}
Therefore, the condition (2.6) is necessary for the origin to be a
center of  equation (2.1).

{\bf Sufficiency.}

Now, we show that  the condition (2.6) is  also sufficient  center
condition.

 Indeed, by Lemma 3.3 and Lemma 3.4,  under the condition (2.6), we have
$$P_{2m+1}=P_1\sum_{k=0}^m\eta_{2k}\bar{P}_1^{2k}=(\bar{P}_1)'\sum_{k=0}^m\eta_{2k}\bar{P}_1^{2k},$$
where $\eta_{2j}\,(j=0,1,2,...,m)$ are constants. By Lemma 2.2,
$r=0$ is a center of equation (2.1).

In summary, the Theorem 2.2 has been proved. $\Box$

\bigskip

 When $n=2m$, because all the calculations and proofs are similar as  Theorem 2.2, we only give  a brief proof of Theorem 2.3 here.

{\bf Proof.} Similar discussion as above,  when $n=2m$, solving
equation
\begin{equation}\sum_{n=1}^{\infty}a_n'(\theta)c^n=P_1(\theta)(\sum_{n=1}^{\infty}a_n(\theta)c^n)^2+P_{2m}(\theta)(\sum_{n=1}^{\infty}a_n(\theta)c^n)^{2m+1}.\label{4.1}\end{equation}
we get
$$
a_k=\tilde{P}_1^{k-1},\,(k=1,2,...,2m)
$$
$$
a_{2m+1+k}=\tilde{P}_1^{2m+k}+\sum_{j=0}^k\lambda_k^j\tilde{P}_1^{k-j}\overline{\tilde{P}_1^jP_{2m}},\,(k=0,1,2,....)
$$
where $\lambda_k^j\, (j=0,1,2,..,k)$ are constants. From
$a_{2m+k+1}(2\pi)=0$ it follows that
$$\int_0^{2\pi}\tilde{P}_1^{j}P_{2m}d\theta=0, (j=0,1,2,...k)$$
and
$$\int_{0}^{2\pi}\bar{P}_1^{2j}P_{2m}d\theta=0,(j=0,1,2,...,m),$$
which means that the  condition (2.7) is necessary for the origin to
be a center of  equation (2.2).

Similar as Lemma 3.3 and Lemma 3.4 we can prove that under condition
(2.7), we have
$$P_{2m}=P_1(\beta_1\bar{P}_1+\beta_3\bar{P}_1^3+...+\beta_{2m-1}\bar{P}_1^{2m-1})=({\bar{P}_1})'(\sum_{j=0}^m\beta_{2j-1}\bar{P}_1^{2j-1}),$$
where $\beta_{2j-1}\,(j=0,1,2,...,m)$ are constants. By Lemma 2.2,
$r=0$ is a center of equation (2.2). $\Box$

\bigskip

{\bf Remark}  By Theorem 2.1, we can derive all the  focal quantity
formulas ($\frac{[n]}{2}+1$) of system (1.5) with arbitrary $n$,
which are more concise and beautiful than the results of calculation
by using computer.

 Taking $n=2,3,4,5,6 $ respectively in (1.6), by
Theorem 2.1, we get the following corollaries.

\begin{col}  $r=0$ is a center of equation
$$\frac{dr}{d\theta}=r(P_1r+P_2r^2)$$
if and only if
$$\int_{0}^{2\pi}\bar{P}_1^{2j}P_2d\theta=0,\,(j=0,1)$$ i.e., the
origin of the cubic system (1.5) $(n=2)$ is a center if and only if
$$p_{20}+p_{02}=0,$$
$$p_{20}(p_{01}^2-p_{10}^2)-p_{11}p_{10}p_{01}=0.$$
\end{col}
This result is the same as Alwash's\cite{5}.

\begin{col}  $r=0$ is a center of equation
$$\frac{dr}{d\theta}=r(P_1r+P_3r^3)$$
if and only if
$$\int_{0}^{2\pi}\bar{P}_1^{2j+1}P_3d\theta=0,\,(j=0,1)$$ i.e., the
origin of the quartic system (1.5) $(n=3)$ is a center if and only
if
$$p_{10}p_{21}-p_{01}p_{12}+3p_{10}p_{03}-3p_{01}p_{30}=0,$$
$$p_{30}p_{01}^3-p_{21}p_{01}^2p_{10}+p_{12}p_{01}p_{10}^2-p_{03}p_{10}^3=0.$$
\end{col}

\begin{col}  $r=0$ is a center of equation
$$\frac{dr}{d\theta}=r(P_1r+P_4r^4)$$
if and only if
$$\int_{0}^{2\pi}\bar{P}_1^{2j}P_4d\theta=0,\,(j=0,1,2)$$ i.e., the
origin of the quintic system (1.5) $(n=4)$ is a center if and only
if
$$3(p_{40}+p_{04})+p_{22}=0,$$
$$(p_{10}^2-p_{01}^2)(p_{04}-p_{40})=p_{10}p_{01}(p_{31}+p_{13}),$$
$$p_{40}p_{01}^4-p_{31}p_{01}^3p_{10}+p_{22}p_{01}^2p_{10}^2-p_{13}p_{01}p_{10}^3+p_{04}p_{10}^4=0.$$
\end{col}

\begin{col}  $r=0$ is a center of equation
$$\frac{dr}{d\theta}=r(P_1r+P_5r^5)$$
if and only if
$$\int_{0}^{2\pi}\bar{P}_1^{2j+1}P_5d\theta=0,\,(j=0,1,2)$$ i.e.,
the origin of the sextic system (1.5) $(n=5)$ is a center if and
only if
$$ p_{01} (5p_{50}+p_{32}+p_{14})=p_{10}(5p_{05}+p_{23}+p_{41}),$$
$$p_{10}^3(p_{23}+10p_{05})-3p_{10}^2p_{01}(2p_{14}+p_{32})+3p_{10}p_{01}^2(p_{23}+2p_{41})-p_{01}^3(p_{32}+10p_{50})=0,$$
$$p_{50}p_{01}^5-p_{41}p_{01}^4p_{10}+p_{32}p_{01}^3p_{10}^2-p_{23}p_{01}^2p_{10}^3+p_{14}p_{01}p_{10}^4-p_{05}p_{10}^5=0.$$
\end{col}

\begin{col}
$r=0$ is a center of equation
$$\frac{dr}{d\theta}=r(P_1r+P_6r^6)$$
if and only if
$$\int_{0}^{2\pi}\bar{P}_1^{2j}P_6d\theta=0,\,(j=0,1,2,3)$$ i.e.,
the origin of the  system (1.5) $(n=6)$ is a center if and only if
$$5(p_{60}+p_{06})+p_{42}+p_{24}=0,$$
$$p_{10}^2(p_{24}-5p_{60}+10p_{06})-p_{10}p_{01}(5p_{51}+5p_{15}+3p_{33})+p_{01}^2(p_{42}-5p_{06}+10p_{60})=0,$$
$$p_{10}^4(p_{24}-3p_{60}+12p_{06})-p_{10}^3p_{01}(3p_{51}+7p_{15}+3p_{33})-12p_{10}^2p_{01}^2(p_{60}+p_{06})$$$$-p_{10}p_{01}^3(3p_{15}+7p_{51}+3p_{33})+p_{01}^4(p_{42}-3p_{06}+12p_{60})=0,$$
$$p_{60}p_{01}^6-p_{51}p_{01}^5p_{10}+p_{42}p_{01}^4p_{10}^2-p_{33}p_{01}^3p_{10}^3+p_{24}p_{01}^2p_{10}^4-p_{15}p_{01}p_{10}^5+p_{06}p_{10}^6=0.$$
\end{col}

{\bf Remark.} By the above corollaries, it can be seen  that the
last focus quantity, i.e., the $[\frac{n}{2}]+1$th  focus quantity,
is very easy to obtain by Theorem 2.1. And the expressions of the
focal values are more succinct and beautiful than the results
calculated by computer.

\section*{ Acknowledgements}

This work  is supported by the  National Natural Science Foundation
of China (61773017,  11571301)  and the National Natural Science
Foundation of Province Jiangsu ( BK20161327).

\end{document}